\documentclass[11pt]{article}
\usepackage{amssymb}
\usepackage{amsmath}
\usepackage{amsthm}

\newtheorem{thm}{Theorem}
\newtheorem{prop}[thm]{Proposition}
\newtheorem{cor}[thm]{Corollary}

\newcommand \al{\alpha}
\newcommand\be{\beta}

\newcommand\ro{\rho}

\newcommand\ph{\varphi}

\newcommand\Om{\Omega}

\def\RR{\mathbb R}

\newcommand\oo{{\infty}}

\renewcommand\o{\circ}
\renewcommand\div{\on{div}}
\newcommand\x{\times}
\newcommand\on{\operatorname}

\newcommand\Ad{\on{Ad}}
\newcommand\ad{\on{ad}}

\newcommand\vol{\on{vol}}

\newcommand\grad{\on{grad}}

\newcommand\Aut{\on{Aut}}
\newcommand\Der{\on{Der}}
\newcommand\Diff{\on{Diff}}

\newcommand\SO{\on{SO}}

\newcommand\g{\mathfrak g}

\newcommand\h{\mathfrak h}

\date{ }
\begin{document}

\title{Geodesics and curvature of semidirect product groups}
\author{Cornelia Vizman \\\it \small West University of Timisoara\\ 
\it\small Bd. V.Parvan 4, 1900-Timisoara, Romania\\
\it \small e-mail: vizman@math.uvt.ro}
\maketitle

\begin{abstract}
Geodesics and curvature of semidirect product groups with right invariant
metrics are determined. In the special case of an isometric semidirect 
product, the curvature is shown to be sum of the curvature of the two groups. 
A series of examples, like the magnetic extension of a group, are 
then considered.
\end{abstract}

\section{Introduction}
\footnotetext{1991 Mathematics Subject Classification: 58D05, 58B20}
\footnotetext{Key words and phrases: semidirect product, geodesic, curvature,
stability, magneto-hydrodynamics}
Important partial differential equations were obtained 
as geodesic equations on diffeomorphism groups with right invariant 
$L^2$ or $H^1$ metrics: 
Euler equation of motion of an incompressible ideal fluid and
the averaged Euler equation
are geodesic equations on the volume preserving diffeomorphism
group [A][MRS],
the Korteweg-de-Vries equation and 
Camassa-Holm shallow water equation are geodesic equations on
the Virasoro-Bott group [OK][M1],
Burger's equation and Camassa-Holm equation are 
geodesic equations on $\Diff(S^1)$ [K],
the superconductivity equation is geodesic equation
on a central extension of the volume preserving diffeomorphism
group [V2].

Semidirect products are applied to study differential equations in physics,
like ideal magneto-hydrodynamics and compressible magneto-hydrodynamics,
heavy top, compressible fluid, elasticity, plasma. In some cases 
one can write the differential equation as a geodesic equation on
semidirect product groups with right invariant metrics:
Kirchhoff equations for a rigid body moving in a fluid are geodesic
equations on the Euclidean group,
the equations of ideal magneto-hydrodynamics are geodesic
equations on the semidirect product of the group of volume preserving
diffeomorphisms and the linear space of divergence free vector fields [ZK],
the equation of passive motion in ideal hydrodynamics is geodesic 
equation on the semidirect product of
the group of volume preserving
diffeomorphisms and the linear space of smooth functions [H].

Arnold suggested an approach to the stability of those differential
equations obtained as geodesic equations on a Lie group with
right invariant metric, by studying the curvature of this weak
Riemannian manifold. The curvature tensor enters the Jacobi equation 
and this, being the linearization of the geodesic equation, controls 
the infinitesimal stability of geodesics. As in the finite-dimensional case, 
one can expect that negative curvature causes exponential instability
of geodesics. 
In this way Arnold [A] showed the instability in most directions 
of Euler equation for ideal flow, Shkoller [S] showed that the
averaged Euler equation is more stable than Euler equation for ideal flow
and Zeitlin and Kambe [ZK] showed that the equations of ideal 
magneto-hydrodynamics are more stable than Euler equation for ideal flow.
There are also results on the stability of the 
Korteweg-de-Vries and Camassa-Holm equations [M1][M2]
and superconductivity equation [V2].

In this paper we determine the geodesics and the curvature in 
the general setting: semidirect product groups with right invariant metrics.
In the special case of an isometric semidirect product $G\ltimes H$,
the curvature is shown to be
the sum of the curvatures of the two groups $G$ and $H$. 
The formulas are applied to several examples like: linear action,
conjugation, $\Diff(M)\ltimes C^\oo(M)$, magnetic extension of a group:
$G\ltimes\g^*$.

\section{Right invariant metrics on Lie groups}

In this paragraph we give expressions for the geodesic equation, Levi-Civita 
covariant derivative and curvature for Lie groups with right invariant 
metrics (see [MR] for a nice presentation of this subject).

Let $G$ be a Lie group with Lie algebra $\g$. Let $\ro_x$ be
the right translation by $x$. Any right invariant bounded Riemannian 
metric on $G$ is determined by its value at the identity $<,>:\g\x\g\to\RR$, 
a positive definite bounded inner product on $\g$.
Let $g:I\to G$ be a smooth curve and $u:I\to\g$ its right logarithmic 
derivative (the velocity field in the right trivialization) 
$u(t)=T\ro_{g(t)^{-1}}.g'(t)$.
In terms of $u$ the geodesic equation for $g$ has the expression (here
$u_t=\frac{du}{dt}$)
\begin{equation*}
u_t=-\ad(u)^\top u,
\end{equation*}
where $\ad(X)^\top$ is the adjoint of $\ad(X)$ with respect to $<,>$, 
if this adjoint does exist. 

The right trivialization induces an isomorphism 
$R:C^\oo(G,\g)\to\mathfrak{X}(G)$ given by $R_X(x)=T\ro_x.X(x)$. In terms of 
this isomorphism, the Levi-Civita covariant derivative is
\begin{equation*}
\nabla^G_X Y=dY.R_X+\frac12\ad(X)^\top Y+\frac12\ad(Y)^\top X-\frac12\ad(X)Y,
\end{equation*}
for $X,Y\in C^\oo(G,\g)$.

The sectional curvature ${\cal K}$ and the Riemannian curvature ${\cal R}$ 
are related by 
\begin{equation*}
{\cal K}(X,Y)=\frac{<{\cal R}(X,Y)Y,X>}{<X,X><Y,Y>-<X,Y>^2},
\end{equation*} 
so the sign of the expression 
$<{\cal R}(X,Y)Y,X>$ determines the sign of the sectional curvature. 
Its expression in the right trivialization is
\begin{align*}
<{\cal R}(X,Y)Y,X>&=\frac14\Vert\ad(X)^\top Y+\ad(Y)^\top X\Vert^2
-\frac34\Vert\ad(X)Y\Vert^2\nonumber\\
&-<\ad(X)^\top X,\ad(Y)^\top Y>\\
&-\frac12<\ad(X)^\top Y,\ad(X)Y>
-\frac12<\ad(Y)^\top X,\ad(Y)X>.\nonumber
\end{align*}

{\bf An example} [A]:
Let $G=\Diff_{vol}(M)$ be the group of volume preserving diffeomorphisms 
of a compact Riemannian manifold $(M,g)$
and $\g=\mathfrak{X}_{vol}(M)$ 
the Lie algebra of divergence free vector fields.
We consider the right invariant metric on $G$ given by the $L^2$ inner 
product $<X,Y>=\int_M g(X,Y)\vol$. Let $\nabla$ be the Levi Civita covariant 
derivative and $R$ the Riemannian curvature tensor. The transpose of $\ad(X)$ 
is $\ad(X)^\top Y=P(\nabla_X Y+(\nabla X)^\top Y)$,
the geodesic equation in terms of the right logarithmic derivative
is Euler equation for ideal flow
\begin{equation*}
u_t=-\nabla_uu-\grad p,\quad\div u=0,
\end{equation*}
the covariant derivative for right invariant vector fields is 
$\nabla^G_XY=P\nabla_XY$ and the curvature 
\begin{align*}
<{\cal R}(X,Y)Y,X>
&=<R(X,Y)Y,X>\nonumber\\
&+<Q\nabla_XX,Q\nabla_YY>-\Vert Q\nabla_XY\Vert^2,
\end{align*}
with $P$ and $Q$ the orthogonal projections 
on the spaces of divergence free respectively gradient vector fields.


\section{Geodesics and curvature of semidirect\\product groups}

Let $G$ and $H$ be Lie groups with right invariant metrics given by positive
definite inner products $<,>$ on their Lie algebras $\g$ and $\h$. 
Let $B$ be an action of $G$ on $H$ by group homomorphisms with
$B:G\x H\to H$ smooth. Then the semidirect product group $G \ltimes H$
is a Lie group with group operation: 
$(g_1,h_1)(g_2,h_2)=(g_1g_2,h_1B(g_1)h_2)$.
We denote by $\be:G\to\Aut\h$ the map defined by $\be(g)=T_eB(g)$ and
$b:\g\to\Der\h$ the differential of $\be$ at the identity. Then $b$ defines
a semidirect product of Lie algebras $\g\ltimes\h$ by
$[(X_1,Y_1),(X_2,Y_2)]=([X_1,X_2],[Y_1,Y_2]+b(X_1)Y_2-b(X_2)Y_1)$ 
and this is the Lie algebra of $G\ltimes H$.
On the semidirect product group we consider the right invariant metric 
given at the identity by the following 
positive definite inner product on its Lie algebra:
$<(X_1,Y_1),(X_2,Y_2)>=<X_1,X_2>+<Y_1,Y_2>$.

We make the following assumptions: 
the transpose of $\ad(X)$, $\ad(Y)$ and $b(X)$ with respect to 
the inner products in $\g$ and $\h$ exist for all $X\in\g$ and $Y\in\h$ 
and there exists a bilinear map $h:\h\x\h\to\g$ defined by the relation
$$
<b(X)Y_1,Y_2>=<h(Y_1,Y_2),X>.
$$
Then the transpose of $\ad(X,Y)$ exists in the semidirect product 
Lie algebra and
$$
\ad(X_1,Y_1)^\top(X_2,Y_2)=(\ad(X_1)^\top X_2-h(Y_1,Y_2),\ad(Y_1)^\top Y_2
+b(X_1)^\top Y_2).
$$

\begin{prop}
With the assumptions above, 
on the semidirect product Lie group $G\ltimes H$ with right invariant metric,
the geodesic equation 
in terms of the right logarithmic derivative $(u,\al):I\to\g\ltimes\h$ is
\begin{gather*}
u_t=-\ad(u)^\top u+h(\al,\al)\nonumber\\
\al_t=-\ad(\al)^\top\al-b(u)^\top\al.
\end{gather*}
\end{prop}

After some computations and using the relation
$$
<h(Y_1,Y_2),\ad(X_1)X_2>=<b(X_1)^\top Y_2,b(X_2)Y_1>
-<b(X_2)^\top Y_2,b(X_1)Y_1>
$$
implied by the Lie algebra homomorphism property of $b$, we get
the following formula for the curvature:

\begin{prop}
Suppose that the map $h$ and 
the transpose of $\ad(X)$, $\ad(Y)$ and $b(X)$ exist. Then
the sign of the sectional curvature in the semidirect product group is given
by the sign of
\begin{multline*}
<\tilde{\cal R}((X_1,Y_1),(X_2,Y_2))(X_2,Y_2),(X_1,Y_1)>=
<{\cal R}^G(X_1,X_2)X_2,X_1>\\
+<{\cal R}^H(Y_1,Y_2)Y_2,Y_1>
+\frac14\Vert h(Y_1,Y_2)+h(Y_2,Y_1)\Vert^2-<h(Y_1,Y_1),h(Y_2,Y_2)>\\
-\frac12<h(Y_1,Y_2)+h(Y_2,Y_1),\ad(X_1)^\top X_2+\ad(X_2)^\top X_1>\\
+<h(Y_1,Y_1),\ad(X_2)^\top X_2>+<h(Y_2,Y_2),\ad(X_1)^\top X_1>\\
+\frac14\Vert b(X_1)^\top Y_2+b(X_2)^\top Y_1\Vert^2
-\frac34\Vert b(X_1)Y_2-b(X_2)Y_1\Vert^2\\
-<b(X_1)^\top Y_1,b(X_2)^\top Y_2>\\
-\frac12<b(X_1)^\top Y_1,b(X_2)Y_2>-\frac12<b(X_2)^\top Y_2,b(X_1)Y_1>\\
+<b(X_1)^\top Y_2,b(X_2)Y_1>+<b(X_2)^\top Y_1,b(X_1)Y_2>\\
-\frac12<b(X_1)^\top Y_2,b(X_1)Y_2>-\frac12<b(X_2)^\top Y_1,b(X_2)Y_1>\\
-<\ad(Y_1)^\top Y_1,b(X_2)^\top Y_2>-<\ad(Y_2)^\top Y_2,b(X_1)^\top Y_1>\\
+\frac12 <\ad(Y_1)^\top Y_2,b(X_1)^\top Y_2
+b(X_2)^\top Y_1-b(X_1)Y_2+b(X_2)Y_1>\\
+\frac12 <\ad(Y_2)^\top Y_1,b(X_1)^\top Y_2
+b(X_2)^\top Y_1-b(X_2)Y_1+b(X_1)Y_2>\\
-\frac12<\ad(Y_1)Y_2,b(X_1)^\top Y_2-b(X_2)^\top Y_1+3b(X_1)Y_2-3b(X_2)Y_1>,
\end{multline*}
where ${\cal R}^G$ and ${\cal R}^H$ denote 
the Riemannian curvature of $G$ and $H$. 
\end{prop}

\begin{cor}
The sign of the sectional curvature of a two-dimensional plane spanned by
$(X_1,0)$ and $(X_2,0)$ is given by the expression
\begin{equation*}
<\tilde{\cal R}((X_1,0),(X_2,0))(X_2,0),(X_1,0)>
=<{\cal R}^G(X_1,X_2)X_2,X_1>.
\end{equation*}
The sign of the sectional curvature of a two-dimensional plane spanned by
$(X,0)$ and $(0,Y)$ is given by the expression
\begin{gather*}
<\tilde{\cal R}((X,0),(0,Y))(0,Y),(X,0)>=
<h(Y,Y),\ad(X)^\top X>\nonumber\\
-\frac14\Vert (b(X)+b(X)^\top)Y\Vert^2
+\frac12\Vert b(X)^\top Y\Vert^2-\frac12\Vert b(X)Y\Vert^2.
\end{gather*}
The sign of the sectional curvature of a two-dimensional plane spanned by
$(0,Y_1)$ and $(0,Y_2)$ is given by
\begin{gather*}
<\tilde{\cal R}((0,Y_1),(0,Y_2))(0,Y_2),(0,Y_1)>
=<{\cal R}^H(Y_1,Y_2)Y_2,Y_1>\nonumber\\
+\frac14\Vert h(Y_1,Y_2)+h(Y_2,Y_1)\Vert^2-<h(Y_1,Y_1),h(Y_2,Y_2)>.
\end{gather*}
\end{cor}

In the special case 
when the action $\be$ of $G$ on $\h$ is by isometries,
the formulas reduce considerably. We call this the {\it isometric case}. 
We have then 
\begin{equation*}
<b(X)Y_1,Y_2>+<Y_1,b(X)Y_2>=0,
\end{equation*}
i.e. $b(X)$ is skew-adjoint and so the bilinear map 
$h:\h\x\h\to\g$ is skew-symmetric. The geodesic equation reduces to
\begin{gather*}
u_t=-\ad(u)^\top u\nonumber\\
\al_t=-\ad(\al)^\top\al+b(u)\al.
\end{gather*}
\begin{prop}
In the isometric case the sign of the sectional curvature of the semidirect
product is given by
\begin{multline*}
<\tilde{\cal R}((X_1,Y_1),(X_2,Y_2))(X_2,Y_2),(X_1,Y_1)>\\
=<{\cal R}^G(X_1,X_2)X_2,X_1>
+<{\cal R}^H(Y_1,Y_2)Y_2,Y_1>.
\end{multline*}
\end{prop}

\begin{proof}
Each $b(X)$ is a derivation of the Lie algebra $\h$ and this implies
$$
<\ad(Y_1)Y_2,b(X)^\top Y_3>=<b(X)Y_2,\ad(Y_1)^\top Y_3>
-<b(X)Y_1,\ad(Y_2)^\top Y_3>.
$$
Using this and proposition 2, the result follows.
\end{proof}

\section{Examples}

\hspace{0.5cm} {\bf 1. Linear action}.
Let $V$ be a vector space with an inner product $<,>$ 
and $B$ a linear action of $G$ on $V$.
The semidirect product Lie group $G\ltimes V$ has the Lie algebra  
$\g\ltimes V$. In this case the geodesic equation becomes
\begin{gather*}
u_t=-\ad(u)^\top u+h(v,v)\\
v_t=-b(u)^\top v.
\end{gather*}
If in addition the action $B$ is by isometries, we get the geodesic equation: 
$u_t=-\ad(u)^\top u, v_t=b(u)v$.

{\bf 2. Conjugation}. Let $B$ be the action of $G$ on $G$ by conjugation.
The Lie algebra of $G\ltimes G$ is $\g\ltimes\g$ with
action $b(X)=\ad(X)$. 
Then the map $h:\g\x\g\to\g$ is 
$h(Y_1,Y_2)=-\ad(Y_1)^\top Y_2$ and the geodesic equation on the semidirect
product written for the right logarithmic derivative $(u,v)$ is
\begin{gather*}
u_t=-\ad(u)^\top u-\ad(v)^\top v\\
v_t=-\ad(v)^\top v-\ad(u)^\top v.
\end{gather*}
If the inner product is Ad-invariant, then $<[X,Y],Z>+<Y,[X,Z]>=0$ and 
the geodesic equation becomes $u_t=0, v_t=[u,v]$.

{\bf 3. Natural action of diffeomorphisms on 
functions}. Let $M$ be a compact Riemannian manifold and $B(\ph)f=f\o\ph^{-1}$
the natural action of the diffeomorphism group on the space of smooth real 
functions. The Lie algebra of $\Diff(M)\ltimes C^\oo(M)$ is
${\mathfrak X}(M)\ltimes C^\oo(M)$ with action $b(X)f=-X(f)$, which is exactly
the Lie algebra of differential operators of the first order on $M$.
We consider the $L^2$ metrics on 
vector fields and on functions, i.e. 
$<X,Y>=\int_Mg(X,Y)\vol$ and $<f_1,f_2>=\int_Mf_1f_2\vol$.
The transpose of $\ad(X)$ in ${\mathfrak X}(M)$ is
$\ad(X)^\top Y=\nabla_XY+(\div X)Y+(\nabla X)^\top Y$ (see [V1]).
Moreover: $b(X)^\top f=X(f)+f\div X$ and
$h(f_1,f_2)=-f_2\grad f_1$.
By proposition 1 we find that the geodesic equation on 
$\Diff(M)\ltimes C^\oo(M)$ is 
\begin{gather*}
u_t=-\nabla_u u-(\div u)u-\frac12\grad g(u,u)-f\grad f\\
f_t=-u(f)-f\div u.
\end{gather*}

If we restrict to the subgroup $\Diff_{vol}(M)\ltimes C^\oo(M)$, we are in the
isometric case. The geodesic equation in this case models 
the passive motion in ideal hydrodynamical flow [H]:
\begin{gather*}
u_t=-\nabla_uu-\grad p\\
f_t=-u(f).
\end{gather*}
From proposition 4 it follows that the sign of the sectional curvature in this 
case is
\begin{multline*}
<\tilde{\cal R}((X_1,f_1),(X_2,f_2))(X_2,f_2),(X_1,f_1)>=
<{\cal R}(X_1,X_2)X_2,X_1>
\end{multline*}
where ${\cal R}$ denotes the curvature of $\Diff_{vol}(M)$ and 
$R$ the curvature of $M$. As a consequence we get the result of [ZK] 
(obtained there in the case of a 2-torus) that
the group manifold $\Diff_{vol}(M)\ltimes C^\oo(M)$ is flat 
in all sections containing 
a direction $(0,f)$.

\section{Magnetic extension of a group}

Let $G$ be a Lie group and $\Ad^*$ the coadjoint action on the dual $\g^*$.
The semidirect product $G\ltimes\g^*$ is called the magnetic extension 
of the group $G$ and it is isomorphic to the cotangent group $T^*G$.

Let $A:\g\to\g^*$ be the operator defined  for $X\in\g$ by 
$$
(A(X),Y)=<X,Y>
$$
for any $Y\in\g$. Here $<,>$ is a fixed inner product on $\g$ and
$(,)$ denotes the pairing between $\g$ and its dual. In the case of
the Lie algebra of divergence free vector fields, $A$ is 
called the inertia operator 
of a fluid and it is invertible on the regular part of the dual space.
The dual space $\g^*$ in this case is naturally isomorphic to the 
quotient space $\Om^1/d\Om^0$ of differential 1-forms modulo exact 1-forms
and $A$ sends a vector field $X$ to the coset of the 1-form $X^\flat$
obtained via the Riemannian metric.

We consider again the general case of a Lie group $G$ with an inner
product on its Lie algebra $\g$, such that the transpose of 
$\ad(X)$ exists for any $X\in\g$. 
We denote by $\g^*_{reg}$ the image of $A$ in $\g^*$, so $A:\g\to\g^*_{reg}$
is an isomorphism. Let $<,>$ be the inner product on $\g^*_{reg}$ 
induced via $A$ by the inner product in $\g$. Next we will restrict on the 
subgroup $G\ltimes\g^*_{reg}$.
The map $h:\g^*_{reg}\x\g^*_{reg}\to\g$ is defined by
$<\ad^*(X)a,b>=<X,h(a,b)>$, so $h(A(Y_1),A(Y_2))=\ad(Y_2)^\top Y_1$.
The coadjoint action on the image of $A$ is 
$\ad^*(X)A(Y)=-A(\ad(X)^\top Y)$, so the transpose of $b(X)=\ad^*(X)$ exists  
and $b(X)^\top A(Y)=-A(\ad(X)Y)$.

\begin{prop}
If the transpose of $\ad(X)$ exists for any $X\in\g$, then
the geodesic equation of $G\ltimes\g^*_{reg}$, 
written for the right logarithmic 
derivative $(u,A(v)):I\to\g\ltimes\g^*_{reg}$, is
\begin{gather*}
u_t=-\ad(u)^\top u+\ad(v)^\top v\nonumber\\
v_t=\ad(u)v.
\end{gather*}
\end{prop}

For $G=\SO(3)$ we obtain Kirchhoff equations for a rigid body
moving in a fluid.
For $G=\Diff_{vol}(M)$ with $L^2$ metric
on its Lie algebra of divergence free vector fields, we obtain as
geodesic equation the equations of ideal magneto-hydrodynamics (see also 
[ZK] [H]):
\begin{gather*}
u_t=-\nabla_uu+\nabla_BB-\grad p\\
B_t=-[u,B].
\end{gather*}

We apply proposition 2 and find that
the sign of the sectional curvature in $G\ltimes\g^*_{reg}$ is given
by the sign of
\begin{multline*}
<\tilde{\cal R}((X_1,A(Y_1)),(X_2,A(Y_2)))(X_2,A(Y_2)),(X_1,A(Y_1))>\\
=<{\cal R}^G(X_1,X_2)X_2,X_1>
+\frac14\Vert \ad(Y_1)^\top Y_2+\ad(Y_2)^\top Y_1\Vert^2\\
-<\ad(X_1)Y_1,\ad(X_2)Y_2>
-<\ad(Y_1)^\top Y_1,\ad(Y_2)^\top Y_2>\\
+\frac12<\ad(Y_1)^\top Y_2+\ad(Y_2)^\top Y_1,
\ad(X_1)^\top X_2+\ad(X_2)^\top X_1>\\
+<\ad(Y_1)^\top Y_1,\ad(X_2)^\top X_2>+<\ad(Y_2)^\top Y_2,\ad(X_1)^\top X_1>\\
+\frac14\Vert \ad(X_1)Y_2+\ad(X_2)Y_1\Vert^2
-\frac34\Vert \ad(X_1)^\top Y_2-\ad(X_2)^\top Y_1\Vert^2\\
-\frac12<\ad(X_1)Y_1,\ad(X_2)^\top Y_2>
-\frac12<\ad(X_2)Y_2,\ad(X_1)^\top Y_1>\\
+<\ad(X_1)Y_2,\ad(X_2)^\top Y_1>+<\ad(X_2)Y_1,\ad(X_1)^\top Y_2>\\
-\frac12<\ad(X_1)Y_2,\ad(X_1)^\top Y_2>
-\frac12<\ad(X_2)Y_1,\ad(X_2)^\top Y_1>.
\end{multline*}

With this formula applied to $G=\Diff_{vol}(M)$ we get information 
on the stability of ideal magneto-hydrodynamics. For example the sign
of the sectional curvature of a two dimensional plane spanned by 
$(0,A(Y_1))$ and $(0,A(Y_2))$ is given by $\frac14\Vert P(\nabla_{Y_1}Y_2
+\nabla_{Y_2}Y_1)\Vert^2-<P\nabla_{Y_1}Y_1,P\nabla_{Y_2}Y_2>$ and the sign
of the sectional curvature of a mixed two dimensional plane is the sign of
\begin{multline*}
<\tilde{\cal R}((X,0),(0,A(Y)))(0,A(Y)),(X,0)>
=<P\nabla_XX,P\nabla_YY>\\
-\frac14\Vert P(\nabla_YX+(\nabla X)^\top Y)\Vert^2
+\frac12\Vert [X,Y]\Vert^2-\frac12\Vert P(\nabla_XY+(\nabla X)^\top Y)\Vert^2.
\end{multline*}
In the case of the 2-torus we can recover the result of [ZK] 
that ideal magneto-hydrodynamics is more stable than ideal hydrodynamics.


\begin{thebibliography}{GBNV}

\bibitem[A]{A}
Arnold, V.,
{\it Sur la geometrie differentielle des groupes de Lie de dimension infinie 
et ses applications a l'hydrodynamique des fluides parfaits}, 
Ann. Inst. Fourier, 16(1966), 319-361.

\bibitem[AK]{AK}
Arnold, V.I.; Khesin, B.A.,
{\it Topological Methods in Hydrodynamics},
Springer-Verlag, 1998.

\bibitem[CE]{CE}
Cheeger, J.; Ebin, D.G.,
{\it Comparison Theorems in Riemannian Geometry},
North-Holland, Amsterdam.

\bibitem[H]{H}
Hattori, Y.,
{\it Ideal magnetohydrodynamics and passive scalar motion as geodesics on 
semidirect product groups},
J. Phys. A, 27(1994), L21-L25.

\bibitem[K]{K}
Kouranbaeva, S.,
{\it The Camassa-Holm equation as a geodesic flow on the diffeomorphism group},
Preprint(1997).


\bibitem[MRS]{MRS}
Marsden, J.E.; Ratiu, T., Shkoller, S.,
{\it The geometry and analysis of the averaged Euler equations and a new 
diffeomorphism group}, Preprint, 1999.

\bibitem[MRW]{MRW}
Marsden, J.; Ratiu, T.; Weinstein, A.,
{\it Semidirect product and reduction in mechanics},
Trans. Am. Math. Soc., 281(1984), 147-177.

\bibitem[MR]{MR}
Michor, P.W.; Ratiu, T.,
{\it On the geometry of the Virasoro-Bott Group},
Journal of Lie Theorie, 8(1998).

\bibitem[M1]{M1}
Misiolek, G.,
{\it A shallow water equation as a geodesic flow on the Bott-Virasoro group},
J. Geom. Phys., 24(1998), 203-208.

\bibitem[M2]{M2}
Misiolek, G.,
{\it Conjugate points in the Bott-Virasoro group and the KdV equation}, 
Proc. Amer. Math. Soc. 125(1997), 935-940.


\bibitem[OK]{OK}
Ovsienko, V.Y.; Khesin, B.A.,
{\it Korteweg-de-Vries superequations as an Euler equation},
Funct. Anal. Appl. 21(1987), 329-331.

\bibitem[S]{S}
Shkoller, S.,
{\it Geometry and curvature of diffeomorphism groups with $H^1$ metric 
and mean hydrodynamics},
Journal of Functional Analysis, 160(1998), 337-365.

\bibitem[V1]{V1}
Vizman, C., 
{\it Geodesics and curvature of diffeomorphism groups},
submitted.

\bibitem[V2]{V2}
Vizman, C.,
{\it Geodesics and stability for central extensions of Lie groups;
the superconductivity equation}, submitted.

\bibitem[ZK]{ZK}
Zeitlin, V., Kambe, T.,
{\it Two-dimensional ideal hydrodynamics and differential geometry},
J. Phys. A, 26(1993), 5025-5031.


\end{thebibliography}
\end{document}